\documentclass[10pt,twocolumn,letterpaper]{article}

\usepackage{amsmath}
\usepackage{amssymb}
\usepackage[style=numeric,giveninits=true,minnames=1,maxnames=6,doi=false,url=false,isbn=false]{biblatex}
\usepackage{url}
\usepackage{booktabs}
\usepackage{titling}
\usepackage[export]{adjustbox}
\usepackage[top=0.75in,bottom=1in,left=0.65in,right=0.65in]{geometry}
\AtBeginBibliography{\small}
\usepackage{graphicx}
\graphicspath{{../figures/}}

\pretitle{\begin{center}\Large\bfseries}
\posttitle{\end{center}}

\renewcommand{\vec}[1]{\mathbf{#1}}
\newcommand{\diff}{\mathrm{d}}
\newcommand{\deriv}[2]{\frac{\diff #1}{\diff #2}}
\newcommand{\pderiv}[2]{\frac{\partial #1}{\partial #2}}

\newcommand{\Dedalus}{{\tt Dedalus}}
\newcommand{\figref}[1]{Fig.\ \ref{#1}}

\usepackage{authblk}

\author[1]{Calum S. Skene}
\author[2,3]{Keaton J. Burns}

\affil[1]{Department of Applied Mathematics, University of Leeds, Leeds LS2 9JT, UK}
\affil[2]{Department of Mathematics, Massachusetts Institute of Technology, Cambridge MA 02139, USA}
\affil[3]{Center for Computational Astrophysics, Flatiron Institute, New York NY 10010, USA}

\begin{document}

\title{Fast automated adjoints for spectral PDE solvers}
\date{\vspace{-2em}}

\twocolumn[
\begin{@twocolumnfalse}

\maketitle

\begin{abstract}
\normalsize
We present a general and automated approach for computing model gradients for PDE solvers built on sparse spectral methods, and implement this capability in the widely used open-source Dedalus framework. 
We apply reverse-mode automatic differentiation to symbolic graph representations of PDEs, efficiently constructing adjoint solvers that retain the speed and memory efficiency of this important class of modern numerical methods. 
This approach enables users to compute gradients and perform optimization for a wide range of time-dependent and nonlinear systems without writing additional code. 
The framework supports a broad class of equations, geometries, and boundary conditions, and runs efficiently in parallel using MPI. 
We demonstrate the flexibility and capabilities of this system using canonical problems from the literature, showing both strong performance and practical utility for a wide variety of inverse problems.
By integrating automatic adjoints into a flexible high-level solver, our approach enables researchers to perform gradient-based optimization and sensitivity analyses in spectral simulations with ease and efficiency.
\end{abstract}

\vspace{0.5cm}

\end{@twocolumnfalse}
]


Numerical modeling has become a ubiquitous tool across every discipline of science and engineering.
However, systematically using forward models to enhance scientific understanding often requires costly optimizations over model inputs and parameters.
\emph{Adjoint state methods} provide an efficient and systematic way to compute model gradients needed for such model-based optimization.
Thus, the ability to compute adjoints forms the foundation of many scientific studies that would otherwise be infeasible.
Indeed, the recent rise of machine learning and artificial intelligence as a transformative technology would not have been possible without adjoint methods, which form the basis of the backpropagation algorithms needed to train models.
Beyond machine learning, adjoint methods have a rich history in scientific computing.
They are an indispensable tool in non-modal stability theory \cite{Trefethen_1993, Schmid_2001}, parametric sensitivity analysis \cite{Luchini_2014, Mensah_2020}, and weather forecasting via four-dimensional data assimilation (4D-Var) \cite{Dimet_1986, Lorenc_1986}.
Adjoint techniques have also been applied in countless optimization settings -- for example, to design aerodynamic structures \cite{Jameson_2022}, mitigate thermoacoustic instabilities \cite{Magri2019}, optimize stellarator designs for nuclear fusion reactors \cite{Paul_2019, Landreman2021}, and solve inverse problems in geophysics \cite{Plessix2006}, among many other applications.

Despite their evident utility, adjoint-based computations remain underutilized in many fields, largely due to the difficulty of computing the adjoint operator, which is often a tedious and error-prone process.
For numerical solvers of partial differential equations (PDEs) (e.g., fluid flow or electromagnetics), there are two main strategies for obtaining adjoints.
One approach is to analytically derive an adjoint system by integrating the governing equations by parts in both space and time; the resulting continuous adjoint equations are then discretized and solved numerically.
Alternatively, the PDE can be discretized first, and then the adjoint of the resulting discrete system can be computed using tools like automatic differentiation (AD).
These are known as the continuous and discrete adjoint approaches, respectively, each with its own advantages and disadvantages (see \cite{giles2000introduction, Ma_2021} for detailed discussions).

Several open-source simulation codes now include built-in adjoint functionality (often termed ``differentiable'' codes), using a mix of continuous and discrete approaches.
This includes {\tt SU\textsuperscript{2}} \cite{Palacios_2013}, {\tt FEniCS} \cite{Fenics}, {\tt Firedrake} \cite{FiredrakeUserManual}, {\tt $\Phi_{\text{Flow}}$} \cite{Holl_2024}, {\tt simsopt} \cite{Landreman2021}, {\tt OpenFOAM} \cite{Openfoam_adjoint_manual}, {\tt Exponax} \cite{koehler2024apebench}, {\tt JaxFluids} \cite{Bezgin2023, Bezgin2025}, {\tt Jax-CFD} \cite{Kochkov_2021, Dresdner_2022}, {\tt Trixi.jl} \cite{ranocha2021adaptive}, and Julia's {\tt SciML} ecosystem \cite{rackauckas2020universal}.
The majority of these PDE solvers focus on specific equation sets and are tied to particular applications or geometries.
Notably, {\tt FEniCS} and {\tt Firedrake} are different in that they are high-level finite-element libraries designed for a wide variety of models.
Using these packages, the governing equations are specified in variational (weak) form via the ``unified form language'' \cite{AlnaesEtal2014}, and then automatically discretized by the library.
Adjoint computations in these libraries are executed automatically using the {\tt dolfin-adjoint} framework \cite{Mitusch2019}.
This high-level approach has led to numerous adjoint-based studies across a wide range of PDE applications.

While finite-element methods excel for complex geometries, global spectral methods are a popular and complementary alternative for problems in simpler geometries, such as turbulent fluid flow in a box or a sphere.
These methods offer exponential accuracy as resolution is increased, and they facilitate efficient elliptic solves for differential-algebraic equations as well as implicit time-stepping for stiff problems \cite{boyd2001chebyshev}.
They are widely used in fundamental research and underpin the world's largest turbulence simulations \cite{yeung2025gpu} and state-of-the-art weather forecasting models \cite{wedi2013fast}.
Beyond classical Fourier-based spectral methods, recent advances have introduced large-scale sparse polynomial spectral algorithms that provide Fourier-like speed and accuracy in simple geometries (Cartesian, cylindrical, and polar coordinates).
These developments have enabled cutting-edge simulations of biophysical, geophysical, and astrophysical flows (e.g., \cite{jackson2023scaling, anders2023photometric, vasil2024solar}).

Here we present a numerical approach to automatically generate discrete adjoints for the entire family of fast global spectral methods, implemented in the open-source \Dedalus{} framework \cite{Burns_2020}.
Like {\tt FEniCS} and {\tt Firedrake}, \Dedalus{} is a high-level library that discretizes and solves user-defined PDEs.
However, \Dedalus{} allows users to specify systems of equations in strong form via a simple interface that is both flexible and accessible to computational scientists across many disciplines.
\Dedalus{} includes many spectral bases and has solver paths for eigenvalue problems (EVPs), linear boundary value problems (LBVPs), nonlinear boundary value problems (NLBVPs), and initial value problems (IVPs).
\Dedalus{} is written in \textit{Python}, with compiled extensions for performance-critical routines, and it automatically handles distributed-memory parallelism via MPI.
These features allow \Dedalus{} simulations to be easily prototyped on a laptop and then ran at scale on a high-performance computing clusters.

In this article, we outline our efficient approach to automating adjoint computations for sparse spectral methods.
We demonstrate its capabilities using \Dedalus{} on several representative adjoint-based optimization problems drawn from a variety of scientific domains.
These examples highlight the efficiency, flexibility, and ease of use of our differentiable extension of \Dedalus{}.
By working with discrete adjoints, our method automatically handles any boundary conditions, constraints, or gauge choices, eliminating the difficulties these pose for continuous-adjoint approaches.
Furthermore, our approach allows sparse spectral solvers to be seamlessly integrated with machine learning frameworks, which inherently require differentiable simulators for training.
Importantly, in all the examples below the gradient computation requires minimal user intervention: only a few additional lines of code are needed to obtain the gradient of any cost functional evaluated on a \Dedalus{} solution.
Overall, our work makes adjoint-based optimization readily accessible to computational scientists using spectral methods across disciplines.


\section*{Spectral adjoints}

\subsection*{Continuous vs.\ discrete}

The fundamental definition of an adjoint is mathematically simple: an adjoint operator $\mathcal{A}^\dagger$ to a linear operator $\mathcal{A}$ is one that satisfies $\langle x, \mathcal{A} y\rangle=\langle \mathcal{A}^\dagger x, y\rangle$ for all $x$ and $y$, where $\langle \cdot, \cdot\rangle$ is a suitably chosen inner product.
Efficiently solving the adjoint of a numerical model's Jacobian yields the model's gradients.
Continuous methods discretize the infinite dimensional adjoint, while discrete methods discretize the forward operator and then take its finite-dimensional adjoint.

The continuous adjoint has the advantage of producing the true adjoint of the governing equations.
However, it generally does not provide the exact gradient of the discrete forward model (even if the forward model is well-resolved), and deriving it -- either manually or automatically -- is difficult for constrained systems or complex boundary conditions.
Implementing a continuous adjoint also requires writing a separate solver in addition to the original model.
The discrete adjoint, on the other hand, yields the correct gradient for the discrete forward model.
This advantage often translates into better convergence in outer-loop optimizations, although it can introduce instabilities when the forward model is under-resolved.

Several studies have manually implemented continuous adjoints for specific PDEs in \Dedalus{} using its symbolic equation-entry interface.
For example, adjoint methods have been used to explore the similarities between minimal seeds and instantons \cite{Lecoanet_2018}, to solve inverse problems \cite{Connor_2024}, and to identify subcritical geodynamo solutions \cite{skene_marcotte_2024}.
In a recent study \cite{Mannix_2024}, continuous and discrete adjoints were compared across several Cartesian problems in \Dedalus{}, highlighting the advantages of the discrete approach -- especially when used with optimization methods that benefit from the discrete adjoint's superior convergence.
In each of these cases, implementing the adjoint solver demanded significant user effort (analytical derivations or manual coding of the discrete adjoint), requiring specialist knowledge -- particularly for bounded domains or high-order time-stepping schemes.

Importantly, the use of AD means discrete adjoints can be generated automatically from the forward model code, contributing to their rising popularity in recent years \cite{Hascoet_2013, Farrell2013, jax2018github, Zhang_2022}.
However, current AD toolchains have limitations in programming language compatibility and hardware support.
Moreover, one often needs to implement custom derivatives or ``chain rules'' for many optimized library functions (such as FFTs and linear algebra routines), since manually providing the Jacobian's action can be more efficient and stable than tracing through the full forward computation.
In fact, for the highly structured computations in many PDE solvers, custom rules might be required for nearly the entire code path, limiting the benefits of AD toolchains under such constraints.

Here we take the approach of leveraging the existing flexible, high-level code structures that enable \Dedalus{} to forward model wide varieties of PDEs.
We compute the adjoints of the various solver types at the outer level, and use the code's internal computational graphs to effectively perform a highly structured, efficient, and platform-independent form of AD to evaluate the adjoints of the low-level operator implementations using their associated chain rules.
Our approach eliminates the barriers associated with manually deriving continuous adjoint systems and produces discrete adjoints without requiring any modifications of the direct code or relying on external AD libraries.
This combination enables \Dedalus{} to automatically compute the discrete adjoint for every available solver type, geometry, and time-stepping scheme in an efficient and platform-independent manner.

\subsection*{Sparse spectral methods}

Sparse spectral methods form finite-dimensional discretizations of PDEs by expanding the unknown solutions in a \emph{trial basis} and projecting the equations against a \emph{test basis}, enforcing the strong-form PDE via a weighted-residual method.
Recent advances in the field have established optimal choices for the test and trial bases to produce maximally sparse discretizations for wide ranges of equations and domains \cite{olver2013fast, vasil2016tensor, vasil2019tensor, lecoanet2019tensor, olver2019sparse, ellison2022gyroscopic}.
Linear operators can be directly discretized into sparse (often narrowly banded) matrices, which are fast to apply or invert.
Nonlinear terms can be efficiently evaluated on a collocation grid using fast transforms (e.g.\ FFTs).
Boundary conditions and other constraints can be applied by modifying the discretized systems, or at the continuous level through the use of Lagrange multipliers and tau terms \cite{burns2022corner}.

Generic solvers for wide classes of BVPs and IVPs can then be readily automated as a series of sparse linear system solves against explicitly computed right-hand sides (RHSs); for instance, this encapsulates both Newton iterations for BVPs and mixed implicit-explicit (IMEX) timestepping for IVPs.
For notational brevity, we will denote such a forward solution as a vector of solution coefficients $\vec{X}$ satisfying the equation $\vec{F}(\vec{X}, \vec{p}) = 0$, where $\vec{p}$ is a vector of problem parameters.
In the Methods section, we provide more detail on how EVPs and IVPs are handled, but we will illustrate the process using this simple notation that more naturally matches the form of BVPs.
The goal of the discrete adjoint program is to efficiently compute gradients of solution properties with respect to the parameters, and to do so by reusing the flexible and efficient computational motifs (sparse direct solvers and fast transforms) that enable the automated forward solution of PDEs with spectral methods.

\subsection*{Efficiently calculating gradients}

\begin{figure*}[ht!]
    \centering
    \includegraphics[width=\textwidth]{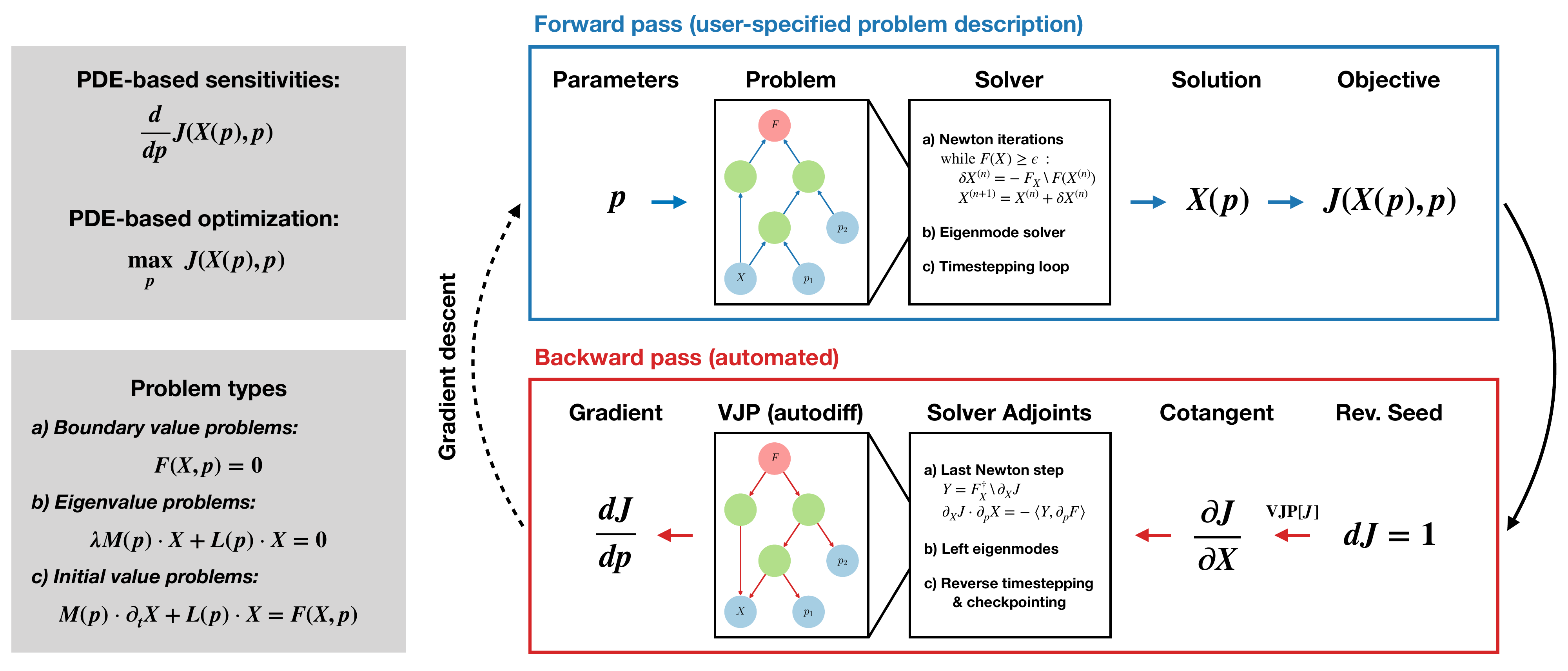}
    \caption{Depiction of the adjoint-looping process to perform efficient PDE-based optimization using automatic discrete adjoints in \Dedalus{}.
    Left: we seek to optimize a functional $J$ of the solution $\vec{X}$ of a PDE (BVP, IVP, or EVP) over parameters $\vec{p}$.
    Top right: \Dedalus{} already enables solving generic PDEs with fast spectral methods with operator graphs representing the PDE (forward pass).
    Bottom right: we have now automated the adjoint (backward pass) by reversing the high-level solver control flow and implementing efficient vector-Jacobian products (VJPs) on the operator graphs, akin to reverse-mode automatic differentiation.
    The result is an optimally fast computation of the discrete model gradient.}
    \label{fig:workflow}
\end{figure*}

Here we summarize how discrete adjoints are used to obtain the derivative with respect to parameters $\vec{p}$ -- which could be forcing terms, initial conditions, or equation parameters -- of a functional applied to the spectral solution of a PDE.
We denote the target functional as $J(\vec{\tilde{X}}, \vec{p})$, where $\vec{\tilde{X}}$ is a generic vector in the finite-dimensional solution space.
We seek to compute the derivative at a point $\vec{\tilde{X}} = \vec{X}(\vec{p})$ which solves the discretized PDE with the specified parameters, i.e.\ $\vec{F}(\vec{X}(\vec{p}), \vec{p}) = 0$.
As a functional purely of the parameters, we define $\mathcal{J}(\vec{p}) = J(\vec{X}(\vec{p}), \vec{p})$.
Directly differentiating this functional by the chain rule yields
\begin{equation}
    \deriv{\mathcal{J}}{\vec{p}} = \pderiv{J}{\vec{p}} + \pderiv{J}{\vec{\tilde{X}}} \cdot \deriv{\vec{X}}{\vec{p}}.
\end{equation}
\noindent where the partial derivatives are evaluated at $(\vec{X}(\vec{p}), \vec{p})$.
The problem in computing this gradient explicitly is the last term, $\diff \vec{X} / \diff \vec{p}$, since directly calculating this Jacobian by e.g.\ finite differences requires solving the forward equation for every element of $\vec{p}$, which becomes prohibitively expensive as the number of parameters increases.

To obtain the gradient with a tractable computational procedure, the adjoint state method approaches this as a constrained optimization problem (for a detailed review, see \cite{Plessix2006}).
The associated Lagrangian is
\begin{equation}
    \mathcal{L}(\vec{\tilde{X}}, \vec{\tilde{Y}}, \vec{p}) = J(\vec{\tilde{X}}, \vec{p}) - \langle \vec{\tilde{Y}}, \vec{F}(\vec{\tilde{X}}, \vec{p})\rangle,
\end{equation}
where $\vec{\tilde{X}}$ and $\vec{\tilde{Y}}$ are generic vectors.
Applying the chain rule directly to $\mathcal{L}$ gives
\begin{equation}
    \deriv{\mathcal{L}}{\vec{p}} = \pderiv{\mathcal{L}}{\vec{p}} + \pderiv{\mathcal{L}}{\vec{\tilde{X}}} \cdot \pderiv{\vec{\tilde{X}}}{\vec{p}} + \pderiv{\mathcal{L}}{\vec{\tilde{Y}}} \cdot \pderiv{\vec{\tilde{Y}}}{\vec{p}}.
\end{equation}
By construction the last term in this expression will be zero when evaluating the gradient at $\vec{\tilde{X}} = \vec{X}(\vec{p})$ because $(\partial \mathcal{L} / \partial \vec{\tilde{Y}}) |_{\vec{X}} = - \vec{F}(\vec{X},\vec{p}) = 0$.
Calculating $\partial \mathcal{L} / \partial \vec{\tilde{X}}$ yields
\begin{equation}
    \pderiv{\mathcal{L}}{\vec{\tilde{X}}} = \pderiv{J}{\vec{\tilde{X}}} - \left(\pderiv{\vec{F}}{\vec{\tilde{X}}}\right)^\dagger \vec{\tilde{Y}},
\end{equation}
where $\dagger$ indicates the Hermitian adjoint.
This term can be set to zero by choosing $\vec{\tilde{Y}} = \vec{Y}(\vec{p})$, where $\vec{Y}(\vec{p})$ satisfies the \emph{adjoint state equation}
\begin{equation}
    \left(\pderiv{\vec{F}}{\vec{\tilde{X}}}\right)^\dagger \vec{Y}(\vec{p}) = \pderiv{J}{\vec{\tilde{X}}},
    \label{equ:adjoint}
\end{equation}
\noindent where again all partials are evaluated at $(\vec{X}(\vec{p}), \vec{p})$.
Hence, by setting $\vec{\tilde{X}} = \vec{X}(\vec{p})$, $\vec{\tilde{Y}} = \vec{Y}(\vec{p})$, and by using the fact that $\mathcal{L}(\vec{X}(\vec{p}), \;\cdot\;, \vec{p}) = \mathcal{J}(\vec{p})$, we obtain
\begin{equation}
    \deriv{\mathcal{J}}{\vec{p}} = \deriv{\mathcal{L}}{\vec{p}} = \pderiv{\mathcal{L}}{\vec{p}} = \pderiv{J}{\vec{p}} - \left\langle \vec{Y}, \pderiv{\vec{F}}{\vec{p}} \right\rangle.
\label{equ:efficient_gradient}
\end{equation}
\eqref{equ:efficient_gradient} is key to the adjoint method and shows that by solving a linear equation for $\vec{Y}(\vec{p})$ we can efficiently obtain the gradient of a functional with respect to parameters, since the dimension of the adjoint state is independent of the number of parameters.

To summarize, in the adjoint state method:
\begin{enumerate}
    \setlength{\itemsep}{0pt}
    \item The forward problem is first solved to determine $\vec{X}(\vec{p})$.
    \item The adjoint problem (\eqref{equ:adjoint}) is then solved for $\vec{Y}(\vec{p})$.
    \item Finally, the right-hand side of \eqref{equ:efficient_gradient} is computed to determine the gradient $\partial \mathcal{J} / \partial \vec{p}$.
\end{enumerate}
The first (forward) step is already automated in \Dedalus{} using optimally sparse spectral methods.
We have now automated the second and third (adjoint) steps, enabling automatic PDE-based optimization through a simple user interface.
The principal technical requirements are:
\begin{itemize}
    \item Implementing fast direct solvers for the Jacobian adjoint systems appearing on the left-hand side of \eqref{equ:adjoint}, based on reusing sparse factorizations from the forward solver.
    \item Adding fast evaluations of the derivative terms on the right-hand side of \eqref{equ:adjoint} and \eqref{equ:efficient_gradient} using matrix-free vector-Jacobian products for arbitrary nonlinear operator trees.
\end{itemize}
The workflow of the adjoint looping process using these tools is illustrated in \figref{fig:workflow}.
More details about the technical implementations are in the Methods section and the Supplementary Information.


\section*{Results}

To demonstrate the simplicity and flexibility of our approach to computing discrete adjoints for sparse spectral solvers, we have created fast implementations using \Dedalus{} of several canonical problems in PDE-based optimization.
We have chosen examples using a variety of solver types, illustrating different scenarios in which gradient information is utilized.
They are also chosen to cover a range of application areas and domain geometries, showcasing the flexibility of our framework and highlighting its suitability for enhancing a wide range of PDE-based modeling studies.
The scripts for these examples are available online\footnote{\url{https://github.com/csskene/dedalus_adjoint_examples}} and include instructions for installing the necessary software.

\subsection*{Parametric sensitivity and numerical continuation}

Since adjoints provide efficient access to gradient information, they are an ideal tool for parametric sensitivity analyses where we seek the direct impact of problem parameters on the PDE solution.
Here we consider the process of computing the neutral stability curve for plane Poiseuille flow \cite{Orszag_1971}, a canonical fluid dynamics problem that describes pressure-driven flow between two flat plates situated at $y=0$ and $y=2$.
In non-dimensional variables, the stability problem can be written as an eigenvalue problem involving the velocity and pressure perturbations $\vec{u}$ and pressure $p$ as
\begin{equation}
\begin{gathered}
    \lambda \vec{u} + \vec{u}_0 \cdot \nabla \vec{u} + \vec{u} \cdot \nabla \vec{u}_0 = -\nabla p + \frac{1}{\textrm{Re}} \nabla^2 \vec{u}, \\
    \nabla \cdot \vec{u} = 0, \\
    \vec{u}(y = \pm 1) = 0,
\end{gathered}
\end{equation}
where $\vec{u}_0 = y(2 - y) \hat{\vec{e}}_x$ is the base flow and $\mathrm{Re}$ is the Reynolds number.
The real part of the eigenvalue $\lambda$ is the growth rate ($\gamma = \Re \, \lambda$), and the imaginary part gives the oscillation frequency.
As the base flow is independent of $x$, this stability problem is separable over streamwise wavenumbers $\alpha$.
For what follows, we consider only 2D perturbations (with no $z$-dependence), and seek to find the sets of parameters ($\alpha$, $\textrm{Re}$) such that the maximum growth rate is zero, indicating neutral stability.

To accelerate the process of computing the neutral curve, we utilize discrete model adjoints in \Dedalus{} in two ways.
First, we guess a point close to the neutral curve and use the EVP solver and its adjoint to compute the largest growth rate $\gamma_\mathrm{max}(\mathrm{Re}, \alpha)$ and its parametric gradient at this point.
Using a Newton solver, we can then efficiently find a point on the neutral curve where $\gamma_\mathrm{max} = 0$.
Second, to construct the guess for a new point on the neutral curve, we use the parametric gradient to find the parametric tangent to the neutral curve, from which a new guess can be estimated.
By repeating this process, we can efficiently trace the neutral curve without difficulties continuing around turning points.

\begin{figure}[ht!]
    \centering
    \includegraphics{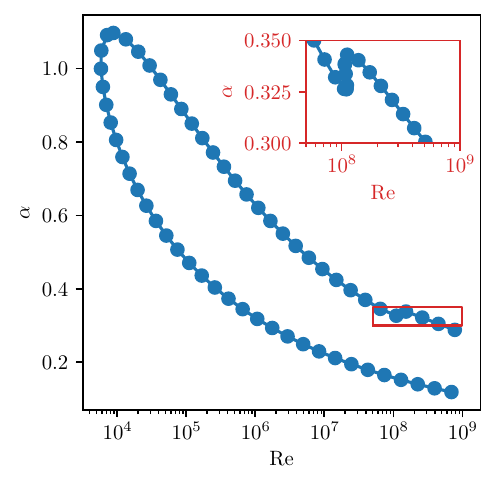}
    \caption{Neutral stability curve for plane Poiseuille flow obtained with adjoint-based eigenvalue sensitivities and numerical continuation.
    Each point on the curve has a growth rate of zero to machine precision.
    The parameters are the horizontal wavenumber $\alpha$ and the base-flow Reynolds number $\mathrm{Re}$.}
    \label{fig:neutral_curve}
\end{figure}

Figure \ref{fig:neutral_curve} shows the result of this procedure.
Each point indicated on this curve has a growth rate of zero to machine precision and was calculated with approximately five eigenvalue solves.
This leads to a very efficient overall calculation without requiring any multidimensional grid searches.
While this equation only contains two parameters, the adjoint method allows quick calculations of eigenvalue sensitivities with respect to potentially many more.
Furthermore, the adjoint implementation reuses the LU decomposition from the forward problem, providing significant computational savings over, e.g., a finite-difference-based approach which would require new forward solves and factorizations.

\subsection*{Nonlinear optimization}

Propagating sensitivities through nonlinear evolutionary equations is possible with \emph{adjoint looping}, where the backward pass corresponds to a reverse-time integration using the time-dependent linearization of the forward model.
This technique allows optimizing, e.g., final-time functionals with respect to equation parameters and initial conditions.
To demonstrate the capabilities of automatic differentiation in \Dedalus{} as a tool for nonlinear optimization, we now consider the problem of optimizing dynamo action in a ball with insulating boundary conditions, following the work of \emph{Chen et al.\ }\cite{Chen_Herreman_Li_Livermore_Luo_Jackson_2018}.
This is an important and challenging example of using nonlinear optimization to assess the stability of fluid flows (see the review \cite{Kerswell_2018} and references therein).

In the astrophysical context, dynamo theory deals with the growth of magnetic fields via the motion of a conducting fluid -- an important physical process underlying magnetic field generation in planets, stars, accretion disks, and galaxies (for a review, see \cite{tobias_2021}).
This problem is governed by the (non-dimensional) induction equation
\begin{equation}
    \begin{gathered}
    \pderiv{\vec{B}}{t} - \eta \nabla^2 \vec{B} = \nabla \times (\vec{u} \times \vec{B}),\\
    \nabla \cdot \vec{B} = 0,
    \end{gathered}
    \label{equ:induction}
\end{equation}
where the fluid velocity is $\vec{u}$ and the magnetic diffusivity is $\eta$.
Although this is a linear equation for $\vec{B}$ when the velocity field is prescribed, as a joint optimization task it is nonlinear (in both parameter and state), and obtaining growing solutions can be a highly nontrivial and subtle task.

Scaling length with the domain radius $R$, velocity with $\Omega R$ for a prescribed vorticity scale $\Omega$, time with $R^2/\eta$, and magnetic field with an arbitrary strength $B$, the induction equation includes a single non-dimensional parameter, the \textit{vorticity}-based magnetic Reynolds number $\textrm{Rm} \equiv R^2 \Omega / \eta$.
In order to ensure the magnetic field remains solenoidal, we formulate the problem using a magnetic vector potential $\vec{A}$, where $\vec{B} = \nabla \times \vec{A}$, as
\begin{equation}
    \begin{gathered}
    \pderiv{\vec{A}}{t} - \nabla^2 \vec{A} + \nabla \psi = \textrm{Rm} \;\vec{u} \times \vec{B},\\
    \nabla \cdot \vec{A} = 0,
    \end{gathered}
    \label{equ:induction_A}
\end{equation}
These equations include a scalar field $\psi$ that enables us to satisfy the Coulomb gauge condition (that the vector potential also be solenoidal).
Since the temporal eigenmodes of $\vec{B}$ correspond to temporal eigenmodes of $\vec{A}$ with the same growth rate, optimizing the growth of $\vec{A}$ beyond initial transients will yield the same optimal velocity field $\vec{u}$ as directly optimizing the growth of $\vec{B}$.

Defining the volume-weighted norm $\|\vec{X}\|^2 = (1/V) \int \vec{X} \cdot \vec{X} \, \textrm{d}V$ and following \cite{Chen_Herreman_Li_Livermore_Luo_Jackson_2018}, we formulate the optimization problem as maximizing
\begin{equation}
    J = \log \|\vec{A}(T)\|^2
\end{equation}
with constraints $\|\vec{A}(0)\| = 1$ and $\|\boldsymbol{\omega}\| = 1$, where $\boldsymbol{\omega} = \nabla \times \vec{u}$ is the vorticity field.
Note that the norm of the vorticity, rather than the velocity, is fixed due to the results of \cite{Proctor_2015}.
We specify vacuum boundary conditions for the magnetic field, which can be written in terms of the spherical harmonic decomposition of $\vec{A}$ as
\begin{equation}
    \pderiv{\vec{A}_{\ell,\sigma}}{r} - \frac{\ell+\sigma}{r} \vec{A}_{\ell,\sigma}  = 0 \quad \textrm{at } r=1,
\end{equation}
where $\ell$ is the spherical harmonic degree and we have used the regularity component index $\sigma \in \{-1,1,0\}$ \cite{vasil2019tensor}.
Since the norm of $\boldsymbol{\omega}$ is constrained, we solve for $\vec{u}$ by solving the auxiliary LBVP
\begin{equation}
    \begin{gathered}
    \nabla \times (\nabla \times \vec{u}) + \nabla \chi = \nabla \times \boldsymbol{\omega},\\
    \nabla \cdot \vec{u} = 0,
    \end{gathered}
\end{equation}
together with no-slip conditions on $\vec{u}$.
This ensures that the velocity field remains divergence-free and satisfies appropriate boundary conditions.
The norm constraints are maintained by using the {\tt Pymanopt} \cite{Townsend_2016} library to perform the optimization on the product manifold $(\vec{A}(0), \boldsymbol{\omega}) \in \Pi = V_{\vec{W}} \times V_{\vec{W}}$, where $V_\vec{W}=\{\vec{x} \;|\; \vec{x}^T \vec{W} \vec{x} = 1\}$ is a generalized Stiefel manifold with a weight matrix $\mathbf{W}$ corresponding to the discretized $L^2$ norm.
Although the IVP is a linear equation for $\vec{A}$, it becomes nonlinear when including $\vec{u}$ as a parameter, and therefore forms a nonlinear optimization problem requiring checkpointing.

\begin{figure*}[ht!]
    \centering
    \includegraphics[height=7cm,valign=T]{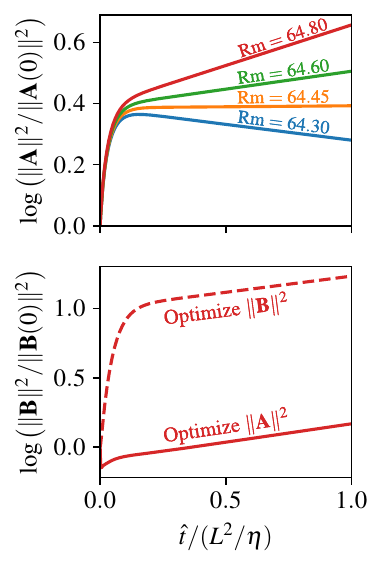}
    \hspace{-0.3cm}
    \includegraphics[width=0.35\textwidth,valign=T]{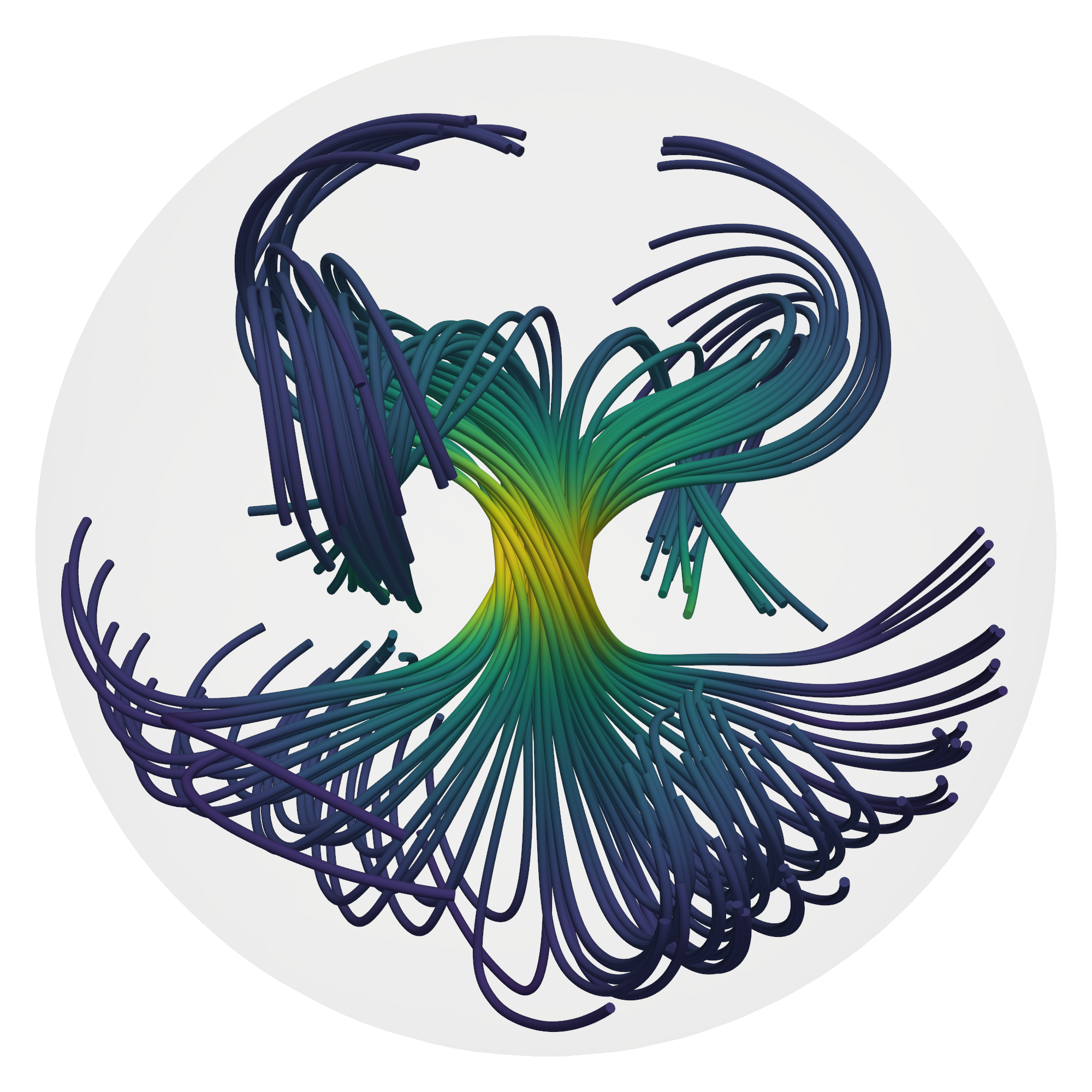}
    \hspace{-0.3cm}
    \includegraphics[width=0.35\textwidth,valign=T]{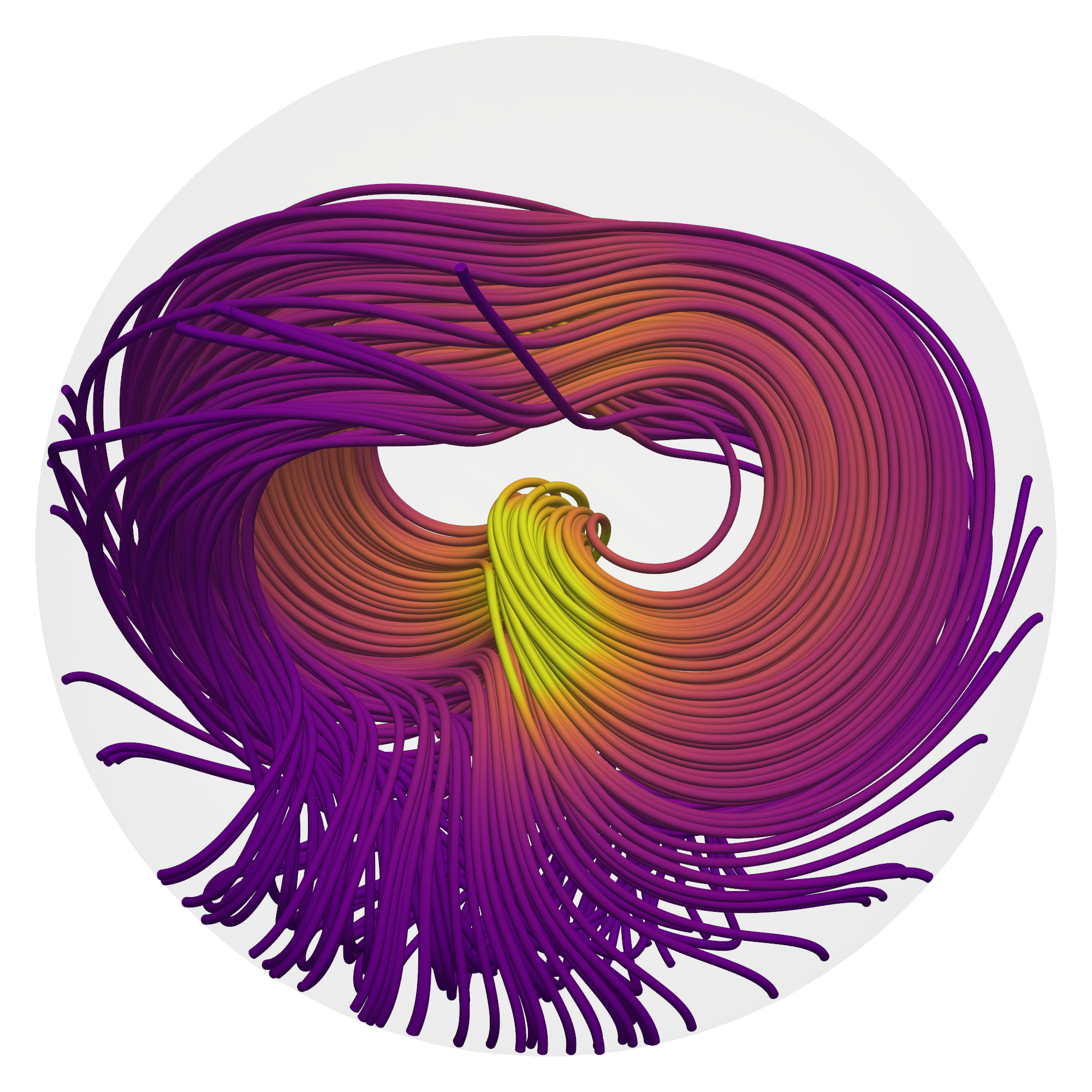}
    \caption{Optimal kinematic dynamo solutions in the ball using nonlinear optimization via adjoint looping.
    Top left: Time series showing the result of optimizing the vector potential norm at different magnetic Reynolds numbers.
    Bottom left: Time series comparing optimization of the magnetic field norm versus the vector potential norm.
    Middle: Streamlines of the optimal velocity field for $\textrm{Rm} = 64.45$.
    Right: Field lines of the optimal magnetic field at $t=1$ for $\textrm{Rm} = 64.45$.}
    \label{fig:optimal_dynamo}
\end{figure*}

The results of the optimization procedure are shown in \figref{fig:optimal_dynamo}.
From the time series, we see that for all magnetic Reynolds numbers there is a short period of transient growth (up to $t \approx 0.1$), followed by an exponentially growing period in which only the mode with the largest growth rate remains.
For $\textrm{Rm} \lesssim 64.45$, the flow is stable and this growth rate is negative, whereas for $\textrm{Rm} \gtrsim 64.45$ the growth rate is positive and exponential growth is observed.
For $\textrm{Rm} \approx 64.45$ the growth rate is approximately zero, indicating that this is the critical magnetic Reynolds number (agreeing with \cite{Chen_Herreman_Li_Livermore_Luo_Jackson_2018}).
The figure also shows the streamlines and field lines of the optimal velocity and magnetic fields, respectively.
We see that the flow is mainly concentrated near the center of the ball, where a large increase in flow velocity occurs with a twisting motion.
Again, this agrees with the results of \cite{Chen_Herreman_Li_Livermore_Luo_Jackson_2018}.

To highlight the flexibility of our method, we also perform the optimization problem of maximizing  
\begin{equation}
    J = \log \| \vec{B}(T) \|^2
\end{equation}
subject to $\|\vec{B}(0)\|_2 = 1$.
This is done by solving the induction equation directly for $\vec{B}$.
\figref{fig:optimal_dynamo} shows the resulting magnetic energy under this optimization compared to the initial optimization for $\vec{A}$.
As expected, both energies show the same growth rates after the initial transient.
However, the amount of transient growth is notably different, with the magnetic-energy-based optimization leading to substantially more absolute growth.
This example emphasizes the flexibility of our automated adjoint capabilities for spectral solvers, allowing the user to easily examine problems from different viewpoints using different equation formulations.

\subsection*{Resolvent analysis}

Resolvent analysis is a powerful modal analysis technique \cite{Taira_2017}, whose applications to turbulence were pioneered by McKeon and Sharma \cite{McKeon_2010}.
In their setup, the flow field is decomposed into a statistically stationary mean flow $\bar{\vec{u}}$ and fluctuations $\vec{u}'$, such that $\vec{u} = \bar{\vec{u}} + \vec{u}'$.
Substituting this into the Navier-Stokes equations and Fourier transforming in time results in the following coupled systems of equations:
\begin{equation}
    \begin{gathered}
    i \omega \hat{\vec{u}}' + \bar{\vec{u}} \cdot \nabla \hat{\vec{u}}' + \hat{\vec{u}}' \cdot \nabla \bar{\vec{u}} + \nabla \hat{p} - \frac{1}{\textrm{Re}} \nabla^2 \hat{\vec{u}}' = \mathcal{F}[\vec{u}' \cdot \nabla \vec{u}'](\omega), \\
    \nabla \cdot \hat{\vec{u}}' = 0,
    \end{gathered}
    \label{equ:fluctuations}
\end{equation}
and
\begin{equation}
    \begin{gathered}
    \bar{\vec{u}} \cdot \nabla \bar{\vec{u}} + \nabla \bar{p} - \frac{1}{\textrm{Re}} \nabla^2 \bar{\vec{u}} = \mathcal{F} , \\
    \nabla \cdot \bar{\vec{u}} = 0,
    \end{gathered}
    \label{equ:mean_flow}
\end{equation}
where we denote the Fourier transform at frequency $\omega$ as $\mathcal{F}[\cdot](\omega)$.
\eqref{equ:fluctuations} show that fluctuations at frequency $\omega$ are driven by triadic interactions, represented by the right-hand side.
In turn, \eqref{equ:mean_flow} show that the zero-frequency component of these triadic interactions sustains the mean flow.

Instead of solving the fully coupled system, mean flow data can be obtained from numerical simulations or experiments in order to directly study the fluctuating components.
In doing so, \eqref{equ:fluctuations} can be reduced to the form:
\begin{equation}
    (i  \omega \vec{M} - \vec{L}) \hat{\vec{u}}' = \hat{\vec{f}},
    \label{equ:resolvent}
\end{equation}
where $\vec{M}$ is a mass matrix, $\vec{L}$ is a linear operator, and the triadic interactions are written as a forcing term $\hat{\vec{f}}$.

Resolvent analysis treats $\hat{\vec{f}}$ as a generic forcing term arising from turbulence and examines \eqref{equ:resolvent} as an input-output system in which a transfer function $\boldsymbol{\mathcal{H}} = (i \omega \vec{M}  - \vec{L})^{-1}$ maps a forcing $\hat{\vec{f}}$ to its driven response $\hat{\vec{u}}'$.
This transfer function $\boldsymbol{\mathcal{H}}$ is known as the resolvent.
Of interest are the leading singular values and vectors of this operator, i.e., $(\sigma_i, \hat{\vec{u}}'_i, \hat{\vec{f}}_i)$ such that $\boldsymbol{\mathcal{H}} \hat{\vec{f}}_i = \sigma_i \hat{\vec{u}}'_i$, and the vectors $\hat{\vec{u}}'_i$ and $\hat{\vec{f}}_i$ are orthonormal under the $L^2$ (energy) norm.
The gains $\sigma_i$ measure the amount of amplification a particular forcing induces in the flow response.
In situations where the largest gain $\sigma_1$ is much greater than the next largest $\sigma_2$, the flow is deemed low-rank, and a generic forcing $\hat{\vec{f}}$ is likely to produce a response that resembles $\hat{\vec{u}}'_1$ due to the orthogonality of the singular basis.

Hence, a resolvent analysis reveals two key insights.
First, it systematically identifies frequencies at which turbulence is likely to manifest, via the magnitude of the gains at each frequency.
Second, for frequencies at which the flow is found to be low-rank, the leading response mode resembles coherent structures expected to be found in fully developed turbulence.
Furthermore, by examining the corresponding forcing mode, one can gain important information relevant for flow control (see, e.g., \cite{Taira_2020, Rolandi_2024}).
While we focus here on resolvent analysis for turbulent flows (about the mean flow), it is worth noting that resolvent analysis for steady flows (about fixed points) is also an important concept in non-modal stability analysis -- predating its use in turbulence \cite{Trefethen_1993, Schmid_2007} -- and serves as the driven counterpart to transient growth analysis.

Using our automatic adjoint routines in \Dedalus{}, we reproduce the resolvent analysis for turbulent pipe flow following \cite{McKeon_2010}.
Since the mean flow is obtained by averaging over $\phi$ and $z$, it depends only on the radial coordinate, allowing us to solve for the forcing and response independently for each azimuthal wavenumber $m$ and streamwise wavenumber $k$.
We use \Dedalus{} to compute the action of both $\boldsymbol{\mathcal{H}}$ and $\boldsymbol{\mathcal{H}}^\dagger$ for given $m$ and $k$ in a cylindrical geometry.
This is done using a linear boundary value problem (LBVP) to apply $\boldsymbol{\mathcal{H}}$ via solving the associated matrix pencil, and applying the adjoint via the vector-Jacobian product of the LBVP.
Together, these routines are used to compute the leading singular vectors via the {\tt SciPy} sparse SVD.
The mean flow is taken from experimental data \cite{McKeon_2004} at $\textrm{Re} = 74345$.

\begin{figure}[ht!]
    \centering
    \includegraphics{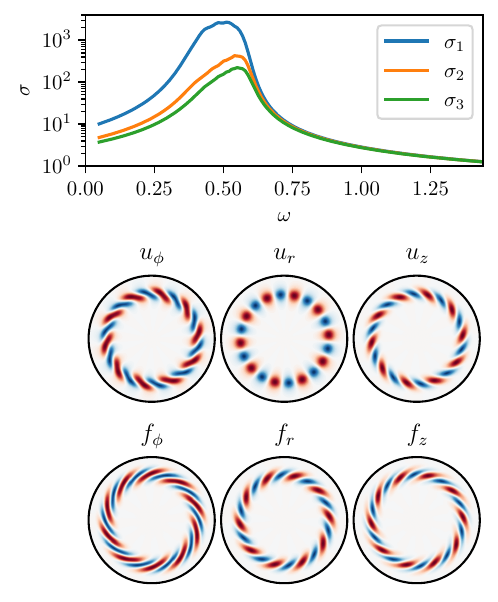}
    \caption{Resolvent analysis for turbulent pipe flow at $\textrm{Re} = 74345$ following \emph{McKeon et al.\ }\cite{McKeon_2004}, for axial wavenumber $k = 1$ and azimuthal wavenumber $m = 10$.
    Top: Leading gains (singular values) as a function of frequency, showing the optimal response at $\omega \approx 0.5$.
    Bottom: Optimal forcing and response modes, showing structures similar to those observed in the experimental fluctuations.}
    \label{fig:resolvent_gains}
\end{figure}

\figref{fig:resolvent_gains} shows the resolvent results for azimuthal wavenumber $m = 10$ and streamwise wavenumber $k = 1$, in excellent agreement with \cite{McKeon_2010}.
The figure illustrates that over a wide range of frequencies, the flow is low-rank, with an order-of-magnitude separation between $\sigma_1$ and $\sigma_2$.
The optimal response near the peak gain at $\omega = 0.5$ is concentrated near the wall and resembles very large-scale motions observed in turbulent pipe experiments \cite{McKeon_2010}.
This example demonstrates the ease of setting up and performing a resolvent analysis using a differentiable spectral code.
Since many turbulent experimental configurations -- such as pipe flow, channel flow, and Taylor-Couette flow -- can be modeled spectrally, there is a clear advantage to complementing these experiments with resolvent analysis.
In particular, we envision that resolvent analyses for fluids not modeled by the incompressible Navier-Stokes equations (e.g., compressible, conducting, non-Newtonian, or viscoelastic fluids) will be made more accessible by these developments, which significantly ease the numerical burden.

\subsection*{Phase reduction analysis}

Phase-reduction analysis is a technique with origins in mathematical biology and neuroscience, used to study phase dynamics and synchronization in oscillating systems \cite{Winfree, Ermentrout2010}.
More recently, it has gained traction in the fluid mechanics community \cite{Taira_Nakao_2018}, where it has been applied to assess the synchronization properties of fluid-structure interactions \cite{Loe_Nakao_Jimbo_Kotani_2021}, thermoacoustic instability \cite{skene_taira_2022}, and to design flow actuation strategies in aeronautical engineering applications \cite{Nair_Taira_Brunton_Brunton_2021, Godavarthi_Kawamura_Taira_2023}.

For a stable limit-cycle solution with period $T$, the notion of phase can be defined as follows.
The phase $\theta$ is first defined for states on the limit cycle as $2\pi t/T$ for $t \in [0, T]$.
In other words, for a state on the limit cycle, the phase variable $\theta \in [0, 2\pi]$ assigns a scalar label to each state, and satisfies the simple ODE $\dot{\theta} = 2\pi/T$ on the cycle.
This definition can be extended to states in the vicinity of the limit cycle by using the fact that the cycle is stable, and therefore any small perturbation returns to it.
Hence, for a state near the limit cycle, we define its phase as the phase of the point on the limit cycle whose trajectory it matches asymptotically as $t \rightarrow \infty$.
Under this definition, the phase dynamics near the limit cycle are governed by the ODE
\begin{equation}
    \dot{\theta} = \frac{2\pi}{T} + \mathbf{z}(\theta) \cdot \vec{h}(\theta),
\end{equation}
where $\vec{h}(\theta)$ represents a small perturbation to the system.
The term $\mathbf{z}(\theta)$ is the phase sensitivity function, whose computation is the central task in a phase-reduction analysis.
By obtaining $\mathbf{z}$, the phase response can be determined independently of the specific form of the perturbation, making this a powerful technique for analyzing the phase and synchronization dynamics of oscillators.

As an example of how adjoints enable spectral methods for phase-reduction analysis, we consider the FitzHugh-Nagumo model \cite{Fitzhugh_1961, Nagumo_1962}:
\begin{equation}
    \begin{gathered}
        \deriv{u}{t} = u - \frac{u^3}{3} - v + I, \\
        \deriv{v}{t} = \epsilon(u + a - b v),
    \end{gathered}
    \label{equ:FHN}
\end{equation}
a simplified version of the Hodgkin-Huxley model \cite{Hodgkin_1952}.
These equations model the firing dynamics of a neuron, with membrane potential $u$ and recovery variable $v$.
The parameters include $\epsilon$ (which sets the timescale separation between $u$ and $v$), the input current $I$, and constants $a$ and $b$ describing the activation dynamics.
For $I \neq 0$, the system exhibits slow-fast repetitive firing for a range of parameter values, yielding a stable limit-cycle solution.

Of interest is the phase response of the dynamics about this limit cycle.
This includes how perturbations can delay or advance the phase of the solution, synchronization of neurons to external stimuli \cite{Kuramoto}, collective synchronization of coupled neurons \cite{Takashi_2004}, and synchronization of oscillators to common white noise \cite{Nakao_2007}.
The key to understanding all of these phenomena lies in computing the phase sensitivity function, which can be efficiently obtained via adjoint methods (see \cite{Ermentrout2010}, for example), as detailed below for the FitzHugh-Nagumo model.
This connection between the adjoint problem and the phase sensitivity function is essential for efficiently conducting a phase-reduction analysis, particularly in PDE systems, as recently demonstrated for fluid flows \cite{Kawamura_2022}.

To compute the phase sensitivity function using \Dedalus{}, we first compute the limit-cycle solution $(u_0, v_0)$.
To do this, we discretize \eqref{equ:FHN} in time using a Fourier basis, and solve it as a nonlinear boundary value problem (NLBVP) to obtain a spectrally accurate periodic solution.
Floquet theory \cite{Floquet1883} then gives that perturbations to this limit cycle $(u', v')$ satisfy the eigenvalue problem
\begin{equation}
    \begin{gathered}
         \deriv{u'}{t} + \lambda u' = u' - u_0^2 u' - v', \\
        \deriv{v'}{t} + \lambda v' = \epsilon(u' - b v'),
    \end{gathered}
    \label{equ:FHN_eig}
\end{equation}
for eigenvalue $\lambda$, where $u'$ and $v'$ are again discretized using a periodic (Fourier) basis in $t$ (the Floquet-Fourier-Hill method \cite{Hill_1886}).
As the system is autonomous, it has a neutral eigenvalue $\lambda = 0$, with corresponding eigenvector $(\dot{u}_0, \dot{v}_0)$ representing phase shifts of the limit cycle.
The phase sensitivity function, which allows the projection of dynamics onto this phase shift, is obtained as the corresponding adjoint eigenvector, normalized such that $\langle \mathbf{z}, (\dot{u}_0, \dot{v}_0)^T \rangle = 1$.
Thus, the phase sensitivity function can be computed in \Dedalus{} simply by solving the adjoint of \eqref{equ:FHN_eig}.

\begin{figure}[ht!]
    \centering
    \includegraphics{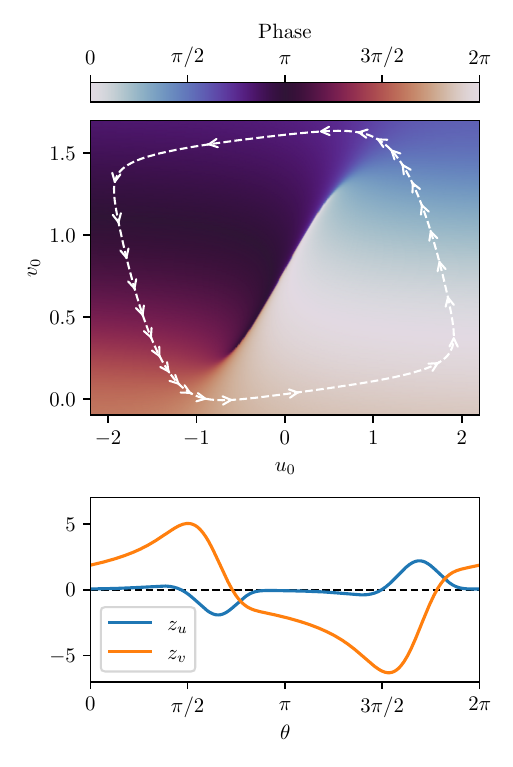}
    \caption{Phase sensitivity analysis of the FitzHugh-Nagumo model via discrete eigenvalue problem adjoints.
    Top: Limit cycle and explicitly computed phase function. Arrows on the limit cycle are equispaced in time, illustrating the slow-fast dynamics of the oscillator.
    Bottom: The phase sensitivity components measuring the gradient of the phase function along the limit cycle.}
    \label{fig:FHN}
\end{figure}

The results are shown in \figref{fig:FHN} for parameters $a = 0.7$, $b = 0.8$, $\epsilon = 0.08$, and $I = 0.8$, demonstrating excellent agreement with \cite{Nakao02042016}.
This approach can be easily adapted to other equation sets, including PDEs, providing a systematic route to phase-reduction analysis for a wide range of oscillators.

\section*{Discussion}

We have presented a framework for efficiently producing adjoint solvers for modern sparse methods and have demonstrated its implementation in the \Dedalus{} library.
By leveraging the structure of spectral solvers, we have developed a hybrid approach in which the adjoint is explicitly provided for core routines and determined via high-level automatic differentiation for others.
By individually handling components such as linear system solves, spectral transforms, timestepping routines, eigenvalue solvers, and Newton solvers, we are able to unobtrusively compute the adjoint of many different PDE solution techniques.

Furthermore, we reuse data structures from the forward code -- such as matrix decompositions and transform plans -- so that the resulting adjoint code is both computationally and memory efficient.
We complement this with automatic differentiation, using it to evaluate vector-Jacobian products of nonlinear terms based on their graph representation in \Dedalus{}.
This enables users to obtain adjoint information for general PDEs directly from the forward model implementation, without requiring additional libraries or restructuring of the original code.

To demonstrate the capabilities of this approach, we have implemented canonical applications of adjoint methods in problems from classical fluid mechanics, dynamo theory, modal analysis of turbulence, and mathematical neuroscience.
These examples span a range of problem types and geometries, showcasing the versatility and advantages of a generic differentiable spectral code.
Specifically, our examples illustrate how such a code can be used for parametric sensitivity analysis, numerical continuation, nonlinear optimization, non-modal stability theory, and phase sensitivity analysis.
These methods are applicable across scientific disciplines and utilize gradient information in diverse ways.
By extending the availability of open-source differentiable solvers to include spectral methods, we provide broad new opportunities for adjoint-based studies in areas where spectral solvers are critical, and do so without requiring custom code development.

While we have focused here on first-order derivative information, a natural next step is the automation of higher-order derivatives \cite{Maddison_2019}.
This would enable techniques such as Newton methods for optimization, as well as uncertainty quantification and Bayesian inference, which can benefit from efficient Hessian-vector products for posterior generation (see, e.g., \cite{Thacker_1989, Isaac_2015}).
Additionally, it is now possible to combine \Dedalus{} models with machine learning libraries by manually linking gradient computations.
Creating a dedicated interface to simplify this integration -- similar to the interface between {\tt Firedrake} and {\tt PyTorch} \cite{Bouziani_2023} -- would enable seamless end-to-end model-based training for a wide range of applications.

Finally, a larger -- but critically important -- future direction is enabling \Dedalus{} solvers to run on graphical and tensor processing units.
Beyond accelerating traditional solvers, this would greatly expand the applicability of \Dedalus{} to machine learning, statistical inference, and nonlinear optimization problems, all of which require solving large numbers of forward and reverse problems.
These exciting research directions build on the strong foundation we have established in rendering \Dedalus{} fully differentiable in a platform-independent fashion.


\section*{Methods}

Solving the adjoint state equation, \eqref{equ:adjoint}, requires adjoint linear-system solves -- e.g., for a sparse matrix $A$ -- and the evaluation of vector-Jacobian products for the right-hand-side (RHS) terms, such as a nonlinear operator $N(\vec{X}, \vec{p})$.
Here, we outline the implementation approach used to efficiently automate these processes in \Dedalus{}.
For details on how different problem types map to this structure, see the Supplementary Information.

\subsection*{Linear system solves}

In each forward iteration, the linear system is solved using a fast direct solver, typically via a sparse LU factorization of $A$.
The corresponding adjoint solve is then obtained by applying the adjoints of the LU factors of $A$.

\subsection*{Nonlinear operator graphs}

Generic nonlinear terms, such as $N$, are represented in \Dedalus{} as directed acyclic graphs (DAGs) of operators acting on the fields and parameters of the PDE.
These graphs are evaluated via depth-first recursion, which triggers spectral transforms as needed to compute intermediate operators.

Forward-mode sensitivity propagation (tangent mode) through an operator graph is implemented as a matrix-free Jacobian-vector product (JVP) evaluated concurrently with the primal graph.
This only requires each operator to provide its discrete derivative.

Reverse-mode sensitivity propagation (adjoint mode, or cotangent propagation) is implemented as a matrix-free vector-Jacobian product (VJP) utilizing the existing DAG operator structure.
During the forward evaluation of the operator tree, a tape records the evaluation order of intermediate operators, forming a topological sort of the DAG.
The VJP is computed by propagating the cotangents backward through this recorded evaluation order and accumulating at the root nodes.
This approach is equivalent to reverse-mode automatic differentiation but is performed in a reliably memory- and compute-efficient manner by leveraging the high-level structure of the forward operator DAG.

\subsection*{Spectral transforms}

During the evaluation of nonlinear operator graphs, spectral transforms are employed to efficiently apply differential operators to field coefficients in their sparse spectral representations, while nonlinearities are evaluated on a collocation grid.
This pseudospectral approach is substantially faster than either pure-collocation or pure-spectral approaches when fast transforms are available.

When evaluating operator VJPs, the adjoint of each spectral transform must therefore be applied.
To automate this process, we have implemented specialized ``adjoint field'' classes that represent cotangent quantities and overload transformation behaviors to automatically apply the appropriate adjoint transforms.

Let $\vec{f}$ denote a field's values on the collocation grid and $\vec{\hat{f}}$ its spectral coefficients.
Let $\vec{T}$ denote the forward transform, such that $\vec{\hat{f}} = \vec{T} \vec{f}$ and $\vec{f} = \vec{T}^{-1} \vec{\hat{f}}$.
Denoting an adjoint field as $\vec{f}^\dagger$, we require that $\langle \vec{f}^\dagger, \vec{f} \rangle = \langle \vec{\hat{f}}^\dagger, \vec{\hat{f}} \rangle$.
Using the above identities, we obtain the relationships:
\begin{equation}
    \vec{\hat{f}}^\dagger = \vec{T}^{-\dagger} \vec{f}^\dagger, \quad \vec{f}^\dagger = \vec{T}^\dagger \vec{\hat{f}}^\dagger.
\end{equation}
This shows that the forward transform of an adjoint field is the adjoint of the original backward transform, and vice versa.

For transforms implemented directly as matrix multiplications, the adjoints are computed using the Hermitian transpose of the matrix.
For fast transforms (e.g., FFTs), the adjoint is implemented using a corresponding fast transform, depending on the transform type.

\section*{Acknowledgements}

We thank Steve Tobias, Jeff Oishi, Geoff Vasil, Daniel Lecoanet, Ben Brown, Rich Kerswell, Patrick Farrell, Colm Caulfield, Andre Souza, Greg Wagner, James Maddison, and Peter Schmid for helpful discussions.

This work was undertaken on ARC4, part of the High Performance Computing facilities at the University of Leeds, UK.
CSS acknowledges partial support from a grant from the Simons Foundation (Grant No. 662962, GF).
CSS would also like to acknowledge support of funding from the European Union Horizon 2020 research and innovation programme (grant agreement no. D5S-DLV-786780).
This work was supported by a grant from the Simons Foundation (SSRFA-6826, K.J.B.).


\newgeometry{left=1in,right=1in,top=0.75in,bottom=1in}

\onecolumn
\appendix{}

\section{Problem types}

Here we give further details of the adjoint approach used for each solver type.

\subsection{Eigenvalue problems}
\texttt{Dedalus} eigenvalue problems take the generalized form
\begin{equation}
\left(\lambda \vec{M}(\vec{p}) + \vec{L}(\vec{p})\right)\vec{X} = 0,
\end{equation}
meaning the adjoint eigenvalue problem is simply
\begin{equation}
\left(\bar{\lambda}\vec{M}^\dagger(\vec{p}) + \vec{L}^\dagger(\vec{p})\right)\vec{Y} = 0.
\end{equation}
For dense solves, the left eigenvectors and eigenvalues can be returned by the eigenvalue library routines wrapped in \texttt{scipy.linalg}.
For sparse solves which find a small number of eigenmodes near a target eigenvalue, a shift-invert method is used, requiring an LU decomposition of a matrix pencil of $\vec{M}$ and $\vec{L}$.
For solving the adjoint eigenvalue problem, this LU decomposition can be reused, leading to significant computational savings for the adjoint solve.

Once the adjoint eigenvalue problem is solved, the sensitivity of the eigenvalue $\lambda$ with respect to parameters is obtained with the equation (see \cite{Luchini_2014}, for example)
\begin{equation}
\pderiv{\lambda}{\vec{p}} = -\frac{\langle \vec{Y},\left(\lambda\pderiv{ \vec{M}}{\vec{p}}+\pderiv{\vec{L}}{\vec{p}}\right)\vec{X}\rangle}{\langle \vec{Y}, \vec{M}\vec{X}\rangle}.
\label{equ:eigenvalue_sens}
\end{equation}
The sensitivity of the eigenvector with respect to parameters is then
\begin{equation}
(\lambda \vec{M}(\vec{p})+\vec{L}(\vec{p}))\pderiv{\vec{X}}{\vec{p}}=-\left( \vec{M}\pderiv{\lambda}{\vec{p}}+\lambda \pderiv{\vec{M}}{\vec{p}} + \pderiv{\vec{L}}{\vec{p}}\right)\vec{X}.
\label{equ:eigenvector_sens}
\end{equation}
This is a singular equation when $\lambda$ is a generalized eigenvalue, but from the Fredholm alternative theorem, it is solvable when the right hand side is orthogonal to the adjoint eigenvector $\vec{Y}$.
This directly gives the eigenvalue sensitivity equation (\ref{equ:eigenvalue_sens}), and therefore there exists a unique solution for the eigenvector sensitivity that is orthogonal to $\vec{X}$.
This sensitivity analysis naturally carries forward to higher order derivatives, as outlined by \cite{Mensah_2020}.

\subsection{Linear boundary value problems}
\texttt{Dedalus} LBVPs take the form
\begin{equation}
\vec{L}\vec{X}=\vec{F}(\vec{p}),
\end{equation}
so the adjoint state equation is simply
\begin{equation}
\vec{L}^\dagger\vec{Y}=\vec{G},
\end{equation}
for a right-hand-side $\vec{G}$ (that depends on the cost functional or upstream solves).
The LU decomposition used for the forward solve can be reused for the adjoint solve. 
A VJP of $\vec{F}$ with cotangents $\vec{Y}$ then propagates sensitivities to the parameters $\vec{p}$ on which $\vec{F}$ depends.

\subsection{Nonlinear boundary value problems}
\texttt{Dedalus} nonlinear boundary value problems can be written as
\begin{equation}
\vec{L}\vec{X}=\vec{F}(\vec{X}, \vec{p}),
\end{equation}
and is solved via Newton-Kantorovich iterations.
The adjoint of this solve is obtained via the linearization of the final solution stage, solving
\begin{equation}
\underbrace{\left(\vec{L}-\pderiv{\vec{F}}{ \vec{X}}\right)^\dagger}_{\vec{H}}\vec{Y} = \vec{G},
\end{equation}
for a right-hand-side $\vec{G}$ (that depends on the cost functional or upstream solves). 
The factorization of $\vec{H}$ from the final Newton iteration is reused for the adjoint state equation.
A VJP of $\vec{F}$ with cotangents $\vec{Y}$ then propagates sensitivities to the parameters $\vec{p}$ and state $\vec{X}$.

\subsection{Initial value problems}

\texttt{Dedalus} IVPs take the form
\begin{equation}
\vec{M}\partial_t\vec{X} + \vec{L}\vec{X} = \vec{F}(\vec{X}, \vec{p}, t),
\end{equation}
which are either solved via a multistep IMEX scheme \cite{Wang_2008} or a Runge-Kutta IMEX scheme \cite{Ascher_1997}.
For both of these timestepping schemes, each step requires computing the explicit terms $\vec{F}$, solving linear systems, and updating the state $\vec{X}$.
To provide the adjoint of this sequence, we have manually implemented the adjoint of the timestepping base classes (for an example of how to derive the adjoint timestepping schemes, see \cite{Mannix_2024}).
An $s$-step multistep IMEX scheme takes the form
\begin{equation}
(a^n_{0}\vec{M}+b^n_{0}\vec{L})\vec{X}_{n}=\sum_{i=1}^{s}\left[c^n_i\vec{F}(\vec{X}_{n-i})-(a^n_i\vec{M}+b^n_i\vec{L})\vec{X}_{n-i}\right],
    \label{equ:imex}
\end{equation}
where the coefficients $(a^n_i, b^n_i, c^n_i)$ depend on the timestep at step $n$, and revert to lower order multistep IMEX schemes for early iterations where not enough steps are computed for higher order schemes. 
The adjoint of this scheme is
\begin{equation}
(a_{0}^n\vec{M}+b_{0}^n\vec{L})^\dagger\vec{Y}_{n}=\sum_{i=1}^s\left(c^{n+i}_i\pderiv{\vec{F}(\vec{X}_{n})}{\vec{X}_{n}}-(a^{n+i}_i\vec{M}+b^{n+i}_i\vec{L})\right)^\dagger\vec{Y}_{n+i}.
 \label{equ:IMEX_adjoint}
\end{equation}
Similarly, a Runge-Kutta IMEX scheme with $s$ stages is
\begin{equation}
(\vec{M}+\textrm{d}t^nH_{i,i}\vec{L})\vec{X}_{n,i}=\vec{M}\vec{X}_{n,0}+\textrm{d}t^n\sum_{j=0}^{i-1}[A_{i,j}\vec{F}_{n,j}-H_{i,j}\vec{L}\vec{X}_{n,j}],
\end{equation}
with $\vec{X}_n=\vec{X}_{n,0}=\vec{X}_{n-1,s}$ and where $\textrm{d}t^n$ is the timestep at iteration $n$. 
This has the adjoint
\begin{equation}
    (\vec{M}+\textrm{d}t^nH_{i,i}\vec{L})^\dagger \vec{Y}_{n,i} =
    \begin{cases}
    \vec{Y}_{n+1,0}, & \textrm{if } i=s,\\
    \sum_{j=i+1}^{s} \textrm{d}t^n\left(A_{j,i}\pderiv{\vec{F}_{n,i}}{\vec{X}_{n,i}}-H_{j,i}\vec{L}\right)^\dagger\vec{Y}_{n,j}, & \textrm{if } 0<i<s,
    \end{cases}
\end{equation}
with
\begin{equation}
    \vec{Y}_{n,0}=\sum_{j=1}^{s} \left(\vec{M}^\dagger+\textrm{d}t^n\left[A_{j,0}\pderiv{ \vec{F}_{n,0}}{x_{n,0}} -H_{j,0}\vec{L}\right]^\dagger\right)\vec{Y}_{n,j},
\end{equation}
and where $\vec{Y}_n=\vec{Y}_{n,0}=\vec{Y}_{n-1,s}$.

In both cases we see that an adjoint step requires VJPs of $\vec{F}$ with cotangents, solving linear systems given via transposes of the forward problem matrices, and updating the cotangents.
The LU decompositions in the direct solve loop are reused where possible, i.e.\ unless they must be rebuilt due to a change in the timestep. 
The VJPs of $\vec{F}$ are performed using the automatic differentiation technique discussed in the manuscript.
As the state at each iteration of the forward solve is required to compute VJPs in the adjoint solve, checkpointing is used using the {checkpoint\_schedules} library \cite{Dolci2024}, which includes many state-of-the-art checkpointing schemes \cite{Stumm_2009, Aupy_2016, Aupy_2017, Pringle_2016, Herrmann_2020, Maddison_2024}.
This allows the user to easily choose between memory and disk based checkpointing.

The adjoint timestepping schemes propagate sensitivities backwards in time from the cotangents defined at the final-time. 
These final-time-condition cotangents stem from VJPs with upstream calculations which can be the cost functional, or other solvers such as LBVPs, NLVPs or other IVPs.
Although it seems that this formulation only provides the sensitivities with respect to each direct state $\vec{X}_n$, we can obtain the sensitivity with respect to constant parameters $\vec{p}$ by adding them to the state $\vec{X}$. 
Similarly, if the cost functional depends on a quantity $\mathcal{J}_I$ integrated over $t$, we can include this in the equation set by including the ODE $\partial_t \mathcal{I} = \mathcal{J}_I$.
This means that the integrated quantity at the final iteration $\mathcal{I}_n$ is obtained with the same accuracy as the forward integration, and its influence on the cotangent calculation is now naturally included through the adjoint timestepping routine.

\section{Example implementation details}
Here we give further details of the implementations for each example.

\subsection{Parametric sensitivity and numerical continuation}
The problem is discretized in the wall-normal direction with Chebyshev polynomials with $256$ modes. 
The streamwise and spanwise directions are parameterized by wavenumbers $\alpha$ and $\beta$. 
A sparse eigenvalue solve is performed with the \texttt{scipy.sparse.eigs} routine which wraps \texttt{ARPACK}.

\subsection{Nonlinear optimization}
The ball is discretized with spin-weighted spherical harmonics in the longitudinal and latitudinal dimensions and one-sided Jacobi polynomials in the radial direction \cite{Vasil_2019}. 
We use spherical harmonics up to degree $\ell=15$, and a radial resolution up to polynomials of maximum degree $31$.
Nonlinear terms are dealiased using a $3/2$-dealiasing rule. A second order semi-implicit multistep BDF scheme is used with fixed timestep $\Delta t=5\times10^{-4}$.
We have checked that increasing the resolution and decreasing the timestep does not significantly affect the results.
The direct solution is provided by solving an LBVP to obtain $\vec{u}$ from $\boldsymbol{\omega}$, followed by an IVP to evolve either $\vec{A}$ or $\vec{B}$.
The adjoint solve requires propagating cotangents back from the cost functional through the IVP and then through the LBVP.

\subsection{Resolvent analysis}
The radial direction in the disk is discretized with Zernike polynomials \cite{Vasil_2016} up to a maximum degree of 127, with the annular and spanwise directions parameterized by wavenumbers $m$ and $k$.
A spectrally accurate weight matrix $\vec{W}$ that discretizes the $L^2$-norm is generated with Zernike quadrature. 
The action of the resolvent matrix $\boldsymbol{\mathcal{H}}$ on a forcing vector $\hat{\vec{f}}$ is obtained via an LBVP.
The optimal basis that maximizes the Rayleigh quotient
\begin{equation}
\sigma=\frac{\|\boldsymbol{\mathcal{H}}\hat{\vec{f}}\|_\vec{W}}{\|\hat{\vec{f}}\|_\vec{W}},
\end{equation}
is found by computing the singular value decomposition (SVD) of $\vec{R}_\vec{M}=\vec{M}^{-1}\boldsymbol{\mathcal{H}}\vec{M}$, where $\vec{M}$ is the Cholesky decomposition of the weight matrix, $\vec{W}=\vec{M}^\dagger\vec{M}$ (trivial here since $\vec{W}$ is diagonal).
Finding the SVD requires matrix-free application routines for both $\boldsymbol{\mathcal{H}}_\vec{M}$ and $\boldsymbol{\mathcal{H}}_\vec{M}^\dagger$, the latter of which requires the adjoint of the LBVP.

\subsection{Phase reduction analysis}
Time is discretized with a complex-Fourier basis with $512$ modes. Dealiasing is performed on a grid twice the size to properly account for cubic nonlinearities.

\section{Verification}
\begin{figure}[h!]
    \centering
    \includegraphics{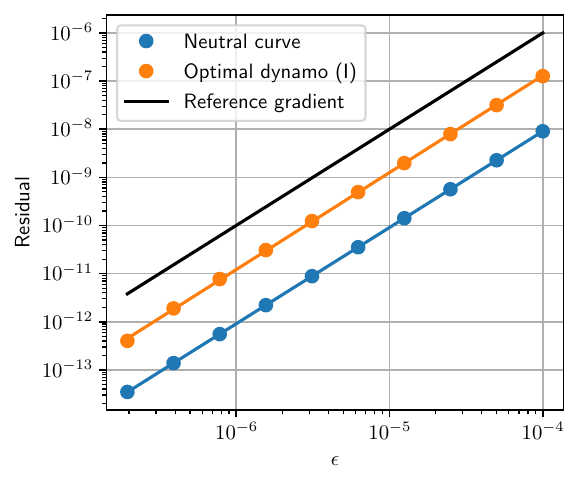}
    \caption{Taylor remainders demonstrating the convergence of finite difference estimates of the adjoint to our automatically computed values. 
    Proper implementation of the discrete adjoint results in second-order convergence (black line).}
    \label{fig:verification}
\end{figure}

In this section we perform verification tests to ensure that our computed gradients are consistent with the direct problem. 
For the plane-Poiseuille and optimal dynamo examples, we use the Taylor-remainder test (see \cite{Farrell2013}, for example).
This test checks that
\begin{equation}
    \mathcal{J}(\vec{p}+\epsilon \vec{p}')=\mathcal{J}(\vec{p}) + \epsilon \langle \vec{Y}, \vec{p}'\rangle + \mathcal{O}(\epsilon^2),
\end{equation}
where $\vec{p}'$ is a perturbation, and $\vec{Y}$ is the computed cotangent. 
This test confirms the correctness of $\vec{Y}$ by showing that
\begin{equation}
|\mathcal{J}(\vec{p}+\epsilon \vec{p}')-\mathcal{J}(\vec{p})-\epsilon \langle \vec{Y}, \vec{p}'\rangle| \approx 0,
\end{equation}
and decays at a second order in $\epsilon$.
This provides a simple and effective test for verifying cotangents and can be applied for linear and nonlinear problems. 
The results of these Taylor-remainder tests are shown in table \ref{tab:verification}, showing excellent accuracy of the computed gradients. 
We also plot the remainders for a range of values of $\epsilon$ in Fig.\ \ref{fig:verification}, showing that the remainders decay at second order over a range of perturbation sizes.

For the linear tests, direct error metrics are shown in Table \ref{tab:verification_linear}, and again indicate that the discrete adjoint is computed correctly.
For the pipe-flow example, we use the fact that $\boldsymbol{\mathcal{H}}_\vec{M}$, is a linear operator to directly verify the adjoint relation
\begin{equation}
\langle \vec{Y}, \boldsymbol{\mathcal{H}}_\vec{M}\vec{X}\rangle =\langle \boldsymbol{\mathcal{H}}_\vec{M}^\dagger\vec{Y}, \vec{X}\rangle.
\end{equation}
For the FitzHugh-Nagumo equation, we check the tangential phase tendency
\begin{equation}
Z_u\deriv{u_0}{t} + Z_v\deriv{v_0}{t} = 1,
\label{eqn:fitzhugh_error}
\end{equation}
for all times. 
This test requires that the phase sensitivity function is the gradient of the phase function $\Phi((u,v)^T)$ evaluated on the limit cycle \cite{Ermentrout2010}. 

\begin{table}[h]
    \centering
    \begin{tabular}{lr}
    Problem &  Taylor-remainder  \\
    \midrule
    Plane-Poiseuille   & 1.999 \\
    Optimal dynamo (I) & 2.016 \\
    Optimal dynamo (II) & 1.993 \\
    \bottomrule
    \end{tabular}
    \caption{Taylor-remainder tests demonstrating the correctness of the discrete adjoint calculations. 
    For the optimal dynamo problem two setups were tested. Setup I: a second-order IMEX BDF scheme with no checkpointing. Setup II: a second-order  Runge-Kutta IMEX scheme using the ``H-revolve'' checkpointing schedule \cite{Herrmann_2020} with 400 checkpoints in RAM and 50 checkpoints on disk.}
    \label{tab:verification}
\end{table}

\begin{table}[h]
    \centering
    \begin{tabular}{lr}
    Problem &  Error  \\
    \midrule
    Pipe-flow   & $1.05\times 10^{-14}$ \\
    FitzHugh-Nagumo & $2.15\times 10^{-13}$ \\
    \bottomrule
    \end{tabular}
    \caption{Direct errors for the linear adjoint examples. For the pipe-flow problem, we report the inner-product test error.
    For the FitzHugh-Nagumo problem, we report the $L^\infty$ error over time of the tangential phase tendency.}
    \label{tab:verification_linear}
\end{table}

\newgeometry{top=0.75in,bottom=1in,left=0.5in,right=0.5in}
\twocolumn
\printbibliography

\end{document}